\documentclass[a4paper,12pt]{article}
\usepackage{latexsym}
\def\beq{\begin{equation}}
\def\eeq{\end{equation}}
\def\bea{\begin{eqnarray}}
\def\eea{\end{eqnarray}}
\def\nn{\nonumber}
\def\bra#1{\left\langle #1\right|}
\def\ket#1{\left| #1\right\rangle}
\def\braket#1#2{\left\langle #1 \left| \right. #2 \right\rangle}
\def\D{$D$-functions }
\def\F{{\cal F}}
\def\Dj#1#2#3{{\cal D}^{\ #1}_{#2,#3}}
\def\Dclass#1#2#3{{\cal D}^{(0)\ #1}_{#2,#3}}
\def\hf{\frac{1}{2}}

\newtheorem{lemma}{Lemma}[section]
\newtheorem{prop}[lemma]{Proposition}
\newtheorem{thm}[lemma]{Theorem}

\newtheorem{cor}[lemma]{Corollary}

\begin{document}
\setlength{\baselineskip}{20pt}
%
%
\begin{center}
\thispagestyle{empty}
\vspace*{3cm}
{\Large \bf Representation Functions for Jordanian Quantum Group 
{\boldmath $SL_h(2)$} and Jacobi Polynomials}

\vspace{2.5cm}
{\large N. Aizawa

\bigskip
\textit{Department of Applied Mathematics \\
Osaka Women's University \\
Sakai, Osaka 590-0035, JAPAN}
}
\end{center}

\vfill
\begin{abstract}
 The explicit expressions of the representation functions ($D$-functions) 
for Jordanian quantum group $ SL_h(2) $ are obtained by combination 
of tensor operator technique and Drinfeld twist. It is shown that the 
$D$-functions can be expressed in terms of Jacobi polynomials as the 
undeformed $D$-functions can. Some of the important properties of the 
$D$-functions for $ SL_h(2) $ such as Winger's product law, recurrence 
relations, RTT type relations are also presented.
\end{abstract}

\newpage
%
%
%
%
\section{Introduction}

  It is known that quantum deformation of Lie group $ GL(2) $ with central 
quantum determinant is classified into two types \cite{ku} : the standard 
deformation $ GL_q(2) $ \cite{frt} and the Jordanian deformation $ GL_h(2) $ 
\cite{demi,zak,ewen}. The representation theory of $ GL_q(2) $ has been 
studied extensively and we know that its contents are quite rich (See, 
for instance, Refs. \cite{bl,cd}). On the other hand, the representation theory 
of $ GL_h(2) $ has not been developed yet. There are some works studying 
differential geometry on quantum $h$-plane and on $ SL_h(2) $ itself \cite{hplane}. 
However, the representation functions for $ GL_h(2) $, the most basic 
ingredient of representation theories, has not been known.  
Recently Chakrabarti and Quesne \cite{cq} showed that the representation functions 
for two-parametric extension of $ GL_h(2) $ \cite{ewen,agha} can be obtained 
from the standard deformed ones via a contraction method and gave explicit form 
of representation functions for some low dimensional cases. 
In Ref.\cite{na1}, the present author shows that the Jordanian deformation 
of symplecton for $ sl(2) $ gives a natural basis for a representation 
of $ SL_h(2) $ and he also gives another basis in terms of quantum $h$-plane. 

  The purpose of the present paper is to obtain explicit formulae for 
$ SL_h(2) $ representation functions using the tensor operator technique 
and to investigate their properties. Representation functions are also called 
Wigner's $D$-functions in physicist's terminology. 
We use both terms and concentrate ourselves to the finite dimensional highest 
weight irreducible representations of $ SL_h(2) $ throughout the present paper. 
In order to make a comparison between \D for $ SL_q(2) $ and 
$ SL_h(2) $, let us recall some known properties of \D for $ SL_q(2) $ \cite{kvs} : 
(a) Wigner's product law \cite{gkk}, (b) recurrence relations \cite{gkk,nom1}, 
(c) orthogonality (d) RTT type relations \cite{nom1}, (e) \D can be written in 
terms of the little $q$-Jacobi polynomials \cite{mmnnu}, 
(f) generating function \cite{nom2}. 
We will show, in this paper, that many of them have counterparts in 
the representation theory of $ SL_h(2) $. Only exception is the generating 
function, it is not presented in this paper. Of course it does not mean that 
the generating function for the \D of $ SL_h(2) $ dose not exist. 

  The plan of this paper is as follows: we present the definitions of 
$ SL_h(2) $ and its dual quantum algebra $ {\cal U}_h(sl(2)) $ in the 
next section. In \S III, before deriving the explicit formulae for the 
representation functions, we discuss general features of them which are 
valid for any kind of deformation of $ SL(2)$ under the assumption that 
the representation theory of the dual quantum algebra has a one-to-one 
correspondence with the undeformed $sl(2).$ Then we shall write down the 
recurrence relations for $ SL_h(2) \ D$-functions. \S IV is a brief review 
of the \D for Lie group $ SL(2) $ (and $ GL(2) $). 
We emphasize that the \D for $ GL(2) $ 
form, in a certain boson realization, irreducible tensor operators of 
the Lie algebra $ gl(2)\oplus gl(2). $ In \S V, a tensor operator technique is used 
to obtain the boson realization of the generators of the Jordanian 
quantum group $ GL_h(2) $, then it is generalized to obtain the 
\D for $ GL_h(2) $. We shall apply the same technique to show that the 
\D for $ SL_h(2) $ can be expressed in terms of Jacobi polynomials. 
This method will be applied to obtain a boson realization for two-parametric 
extension of the Jordanian deformation of $ GL(2) $ in \S VI. \S VII is concluding 
remarks. 

%
%
%
%
\setcounter{equation}{0}
\section{{\boldmath $ SL_h(2) $} and its Dual}

  The Jordanian quantum group $ GL_h(2) $ is generated by four 
elements $ x, y, u $ and $ v $ subject to the relations \cite{demi,zak,ewen}
\bea
 & & [v,\; x] = hv^2, \qquad [u,\; x] = h(D-x^2), 
     \nn \\
 & & [v,\; y] = h v^2, \qquad [u,\; y] = h(D-y^2),
     \label{SLh2} \\
 & & [x,\; y] = h(xv - yv), \qquad [v,\; u] = h(xv + vy),
     \nn
\eea
where $D = xy -uv -hxv$ is the quantum determinant generating 
the center of $ GL_h(2) $. This is a Hopf algebra and Hopf algebra 
mappings have a similar form as $ GL_q(2) $. However, explicit form of the 
mappings is not necessary in the following discussion. By setting 
$ D = 1 $, we obtain $ SL_h(2) $ from $ GL_h(2). $ 

  The quantum algebra dual to $ GL_h(2) $ is denoted by $ {\cal U}_h(gl(2)) $, 
and defined by the same commutation relations as the Lie algebra $ gl(2) $
\beq
   [J_0,\; J_{\pm}] = \pm 2 J_{\pm}, \qquad 
   [J_+,\; J_-] = J_0, \qquad
   [Z,\; \bullet \; ] = 0.              \label{gl2comm}
\eeq
However, their Hopf algebra mappings are modified via twisting \cite{dri} by the 
invertible element $ \F \in {\cal U}_h(gl(2))^{\otimes 2} $ \cite{ogi}
\beq
 \F = \exp\left(-\frac{1}{2} J_0 \otimes \sigma \right), 
 \qquad
 \sigma = -\ln(1 -2hJ_+).                        \label{sigma}
\eeq
The coproduct $ \Delta $, counit $ \epsilon $ and antipode $S$ for 
$ {\cal U}_h(gl(2)) $ are obtained from those for $ gl(2) $ by 
\beq
   \Delta = \F \Delta_0 \F^{-1}, \qquad 
   \epsilon = \epsilon_0, \qquad 
   S = \mu S_0 \mu^{-1},                          \label{hopf}
\eeq
where the mappings with suffix $0$ stand for the Hopf algebra mappings for 
$ gl(2) $. The elements $\mu$  and $\mu^{-1} $ are defined, using the product $m$ for 
$ gl(2) $,  by
\beq
 \mu = m(id \otimes S_0)(\F), \qquad \mu^{-1} = m(S_0 \otimes id)(\F^{-1}).
 \label{mu}
\eeq 
The twist element $ \F $ is not depend on the central element $Z$ so that 
the Hopf algebra mappings for  $Z$ remain undeformed. 
Therefore the Jordanian quantum algebra obtained by the twist element (\ref{sigma}) 
has the decomposition 
$ {\cal U}_h(gl(2)) = {\cal U}_h(sl(2)) \oplus u(1) $. 
The Jordanian quantum algebra $ {\cal U}_h(gl(2)) $ is a triangular Hopf 
algebra whose universal $R$-matrix is given by 
$ {\cal R} = \F_{12} \F^{-1} $. 

  It is obvious, from the commutation relation (\ref{gl2comm}), that 
$ {\cal U}_h(gl(2)) $ and $ gl(2) $ have the same finite dimensional 
highest weight irreducible representations. 
Furthermore we can easily see 
that tensor product of two irreducible representations (irreps) is completely reducible 
and decomposed into irreps in the same way as $ gl(2), $ 
since the Clebsch-Gordan coefficients (CGC)
for $ {\cal U}_h(gl(2)) $ are product of the ones for $ gl(2) $ and matrix elements 
of the twist element $ \F $. For the $ {\cal U}_h(sl(2)) $ sector, this is 
carried out in Ref.\cite{na1}. The CGC for $ {\cal U}_h(sl(2)) $ in 
another basis are discussed in Ref.\cite{vdj}

  Let $ \Delta, \epsilon $ be the coproduct and counit for $ GL_h(2) $, respectively. 
We use the same notations for the Hopf algebra mappings of both $ GL_h(2) $ 
and $ {\cal U}_h(gl(2)) $, however, this may not cause serious confusion. 
A vector space (representation space) $ V $ is called right  
$ GL_h(2) $ comodule, if there exist a map 
$ \rho : V \rightarrow V \otimes GL_h(2)  $ such that the following 
relations are satisfied
\beq
  (\rho \otimes id)\circ \rho = (id_V \otimes \Delta) \circ \rho, 
  \qquad
  (id_V \otimes \epsilon) \circ \rho = id_V,
  \label{comodule}
\eeq
where $ id_V $ stands for the identity map in $V$. The left comodule is 
defined in the similar manner. 
Using the bases $ \{ e_i\ | \ i = 1,2, \cdots, n \} $ of $V$, the map 
$ \rho $ is written as 
\beq
 \rho(e_i) = \sum_j\; e_j \otimes {\cal D}_{ji},     \label{rightcom}
\eeq
it follows that the relation (\ref{comodule}) are rewritten as
\beq
   \Delta({\cal D}_{ij}) = \sum_k {\cal D}_{ik} \otimes {\cal D}_{kj}, 
   \qquad
   \epsilon({\cal D}_{ij}) = \delta_{ij}.
   \label{dfunc}
\eeq
We call $ {\cal D}_{ij} \in GL_h(2) $ satisfying 
(\ref{rightcom}) and (\ref{dfunc}) the $D$-function for $ GL_h(2) $. 

%
%
%
%
\setcounter{equation}{0}
\section{Properties of \D}
\subsection*{A. Wigner's Product Law and RTT Type Relations}

  Before deriving the explicit formulae for $SL_h(2) $ \D, one can discuss 
some important properties of \D such as Wigner's product law, recurrence 
relations, RTT type relations and so on, using the definition of 
universal $T$-matrix \cite{dri0,fg}. The explicit expression of the universal 
$T$-matrix is not necessary. The universal $T$-matrix for the standard 
deformation of $ GL(2) $ is given in Ref.\cite{fg}, while it is not 
known for the Jordanian deformation of $ GL(2) $. 

  The discussion in this subsection is quite general. We shall present it 
so as to be applicable to any kind of deformation of $ SL(2) $ (standard, 
Jordanian, two-parametric extension, anything else if any). Then we will 
write down the results explicitly for the Jordanian deformation of $ SL(2) $ 
in the next subsection. It will also be 
seen that the discussion is easily extended to other groups.

  Let $ {\cal G} $ and  ${\bf g}$ be deformation of Lie group $SL(2) $ and 
Lie algebra $ sl(2) $, respectively. 
The duality between $ {\cal G} $ and $ {\bf g} $ are expressed, 
by choosing suitable bases, in terms of the universal $T$-matrix \cite{fg}. 
Let $ x^{\alpha} $ and $ X_{\alpha} $ be elements of a basis of $ {\cal G} $ 
and $ {\bf g} $, respectively. They are chosen as follows: the product is 
given by
\beq
  x^{\alpha} x^{\beta} = \sum_{\gamma}\; h_{\gamma}^{\alpha, \beta} x^{\gamma},
  \qquad
  X_{\alpha} X_{\beta} = \sum_{\gamma}\; f_{\alpha,\beta}^{\gamma} X_{\gamma},
  \label{dualprod}
\eeq
the coproduct is given by
\beq
 \Delta(x^{\alpha}) = \sum_{\beta,\gamma}\; f_{\beta,\gamma}^{\alpha} 
 x^{\beta} \otimes x^{\gamma},
 \qquad
 \Delta(X_{\alpha}) = \sum_{\beta,\gamma}\; h_{\alpha}^{\beta,\gamma} 
 X_{\beta} \otimes X_{\gamma}.
 \label{dualcopro}
\eeq
Then the universal $T$-matrix $ {\cal T} $ is defined by
\beq
  {\cal T} = \sum_{\alpha}\; x^{\alpha} \otimes X_{\alpha}.
  \label{uniT}
\eeq

  We assume that the deformed algebra $ {\bf g} $ has the same 
finite dimensional highest weight irreps as $ sl(2) $, that is, 
(1) each irrep is classified by the spin $j$ and a irrep basis 
$ \ket{jm} $ is specified by $j$ and the magnetic quantum number $m$, 
(2) tensor product of irreps $j_1$ and $j_2$ is completely reducible
\[
  j_1 \otimes j_2 = j_1+j_2 \oplus j_1+j_2-1 \oplus \cdots \oplus |j_1-j_2|.
\]
We further assume that vectors $ \ket{jm} $ are complete and orthonormal. 
Then the \D for 
$ {\cal G} $ is obtained by
\beq
 {\cal D}^{\ j}_{m',m} = \bra{jm'} {\cal T} \ket{jm} = 
 \sum_{\alpha}\; x^{\alpha} \bra{jm'}X_{\alpha} \ket{jm}.
 \label{DandT}
\eeq
For the standard two-parametric deformation of $ GL(2) $, 
the RHS of (\ref{DandT}) was computed and it was shown that 
(\ref{DandT}) coinceided with the \D obtained by another method \cite{jvdj}. 
In our case, we show that the \D (\ref{DandT}) satisfy (\ref{dfunc}) 
by making use 
of the relations (\ref{dualprod}) and (\ref{dualcopro}). 
The coproduct of $ \Dj{j}{m'}{m} $ is computed as
\bea
 \Delta(\Dj{j}{m'}{m}) &=& \sum_{\alpha}\; \Delta(x^{\alpha}) \bra{jm'} X_{\alpha} \ket{jm} 
 = \sum_{\beta,\gamma}\; x^{\beta}\otimes x^{\gamma} \bra{jm'} X_{\beta} X_{\gamma} \ket{jm}
 \nn \\
 &=& \sum_{\beta,\gamma,k}\; x^{\beta} \otimes x^{\gamma} \bra{jm'} X_{\beta} 
 \ket{jk} \bra{jk} X_{\gamma} \ket{jm} 
 = \sum_k\; \Dj{j}{m'}{k} \otimes \Dj{j}{k}{m}. \nn
\eea
To compute the counit for $ \Dj{j}{m'}{m} $, we use the identiy obtained from 
the definition of counit
\beq
  \sum_{\beta,\gamma}\; f^{\alpha}_{\beta,\gamma} \; \epsilon(x^{\beta}) x^{\gamma} 
  = x^{\alpha}.
  \label{counitid}
\eeq
Using this relation, the universal $T$-matrix is rewritten as
\bea
 {\cal T} &=& \sum_{\alpha}\; x^{\alpha} \otimes X_{\alpha} = 
 \sum\; f^{\alpha}_{\beta,\gamma} \; \epsilon(x^{\beta}) x^{\gamma} \otimes X_{\alpha} 
 = \sum\; \epsilon(x^{\beta}) x^{\gamma} \otimes X_{\beta} X_{\gamma}
 \nn \\
 &=& ( \sum_{\beta}\; \epsilon(x^{\beta}) \otimes X_{\beta} )\; {\cal T}.\nn
\eea
It follows that
\beq
   ( \sum_{\beta}\; \epsilon(x^{\beta}) \otimes X_{\beta} ) 
   = (\epsilon \otimes id) ({\cal T}) = 1. 
   \label{epsionT}
\eeq
Therefore the counit for \D is
\[
 \epsilon(\Dj{j}{m'}{m}) = \bra{jm'} (\epsilon \otimes id) ({\cal T}) \ket{jm} 
 = \braket{jm'}{jm} = \delta_{m',m}.
\]

 We first show that the \D (\ref{DandT}) satisfy the analogous relations to 
 Wigner's product law. 
Let us denote the CGC for $ {\bf g} $ by 
$ \Omega^{j_1,j_2,j}_{m_1,m_2,m} $, $i.e.$, 
\beq
 \ket{(j_1j_2)jm} = \sum_{m_1,m_2}\; \Omega^{j_1,j_2,j}_{m_1,m_2,m} 
 \ket{j_1 m_1} \otimes \ket{j_2 m_2}. 
 \label{cgcg}
\eeq
We write the inverse of the above relation as follows:
\beq
  \ket{j_1 m_1} \otimes \ket{j_2 m_2} = \sum_{j,m}\; 
  \mho^{j_1,j_2,j}_{m_1,m_2,m} \ket{(j_1j_2)jm}.
  \label{mhog}
\eeq
Then an analogue of Wigner's product law reads
\begin{thm} 
The \D for $ {\cal G} $ satisfy the relation
\beq
  \delta_{j,j'} \Dj{j}{m'}{m} = \sum_{k_1,k_2,m_1,m_2}\; 
  \mho^{j_1,j_2,j'}_{k_1,k_2,m'} \; \Omega^{j_1,j_2,j}_{m_1,m_2,m} 
  \Dj{j_1}{k_1}{m_1} \Dj{j_2}{k_2}{m_2}. 
  \label{Wigpro}
\eeq
\label{thmWigpro}
\end{thm}
$ Proof $ : Because of the relations (\ref{dualprod}) and (\ref{dualcopro}), 
one can show that
\beq
 (id \otimes \Delta)({\cal T}) = \sum_{\alpha, \beta}\; 
 x^{\alpha} x^{\beta} \otimes X_{\alpha} \otimes X_{\beta}, 
 \quad
 (\Delta \otimes id)({\cal T}) = \sum_{\alpha, \beta}\; 
 x^{\alpha} \otimes x^{\beta} \otimes X_{\alpha} X_{\beta}.
 \label{idT}
\eeq
It follows that 
\beq
  (id \otimes  \Delta)({\cal T}) \ket{(j_1j_2)jm} = \sum_{\alpha,\beta,m_1,m_2}\; 
  \Omega^{j_1,j_2,j}_{m_1,m_2,m} \; x^{\alpha} x^{\beta} \otimes 
  X_{\alpha} \ket{j_1 m_1} \otimes X_{\beta} \ket{j_2 m_2}.
  \label{prodT}
\eeq
The LHS of (\ref{prodT}) is rewritten as
\[
  \sum_{m'}\; \Dj{j}{m'}{m} \otimes \ket{(j_1j_2)jm'} = 
  \sum_{m',k_1,k_2}\; \Omega^{j_1,j_2,j}_{k_1,k_2,m'} \Dj{j}{m'}{m} \otimes 
  \ket{j_1 k_1} \otimes \ket{j_2 k_2}. 
\]
The RHS of (\ref{prodT}) is rewritten as
\bea
 & & \sum \Omega^{j_1,j_2,j}_{m_1,m_2,m} \; \bra{j_1 k_1} X_{\alpha} \ket{j_1 m_1} 
 \; \bra{j_2 k_2} X_{\beta} \ket{j_2 m_2} 
 x^{\alpha} x^{\beta} \otimes \ket{j_1 k_1} \otimes \ket{j_2 k_2}
 \nn \\
 &=& \sum \Omega^{j_1,j_2,j}_{m_1,m_2,m} \; \Dj{j_1}{k_1}{m_1} 
 \Dj{j_2}{k_2}{m_2} \otimes \ket{j_1 k_1} \otimes \ket{j_2 k_2}.
 \nn
\eea
Thus we obtain
\beq
  \sum_{m'}\; \Omega^{j_1,j_2,j}_{k_1,k_2,m'} \Dj{j}{m'}{m} = 
  \sum_{m_1,m_2}\; \Omega^{j_1,j_2,j}_{m_1,m_2,m} \; \Dj{j_1}{k_1}{m_1} 
 \Dj{j_2}{k_2}{m_2}. 
 \label{prod1}
\eeq
Using the orthogonality of $ \Omega^{j_1,j_2,j}_{m_1,m_2,m} $ and 
$ \mho^{j_1,j_2,j}_{m_1,m_2,m}$, the theorem is proved. \hfill $ \Box $
\begin{cor} 
The \D also satisfy the following relations
\bea
  & & \sum_{m'}\; \Omega^{j_1,j_2,j}_{k_1,k_2,m'} \Dj{j}{m'}{m} = 
  \sum_{m_1,m_2}\; \Omega^{j_1,j_2,j}_{m_1,m_2,m} \; \Dj{j_1}{k_1}{m_1} 
  \Dj{j_2}{k_2}{m_2}, 
  \label{rel1}
  \\
  & & \sum_m\; \mho^{j_1,j_2,j}_{m_1,m_2,m} \Dj{j}{m'}{m} = 
  \sum_{k_1,k_2}\; \mho^{j_1,j_2,j}_{k_1,k_2,m'} \Dj{j_1}{k_1}{m_1} 
  \Dj{j_2}{k_2}{m_2},
  \label{rel2}
  \\
  & & \Dj{j_1}{k_1}{m_1} \Dj{j_2}{k_2}{m_2} = \sum_{j,m,m'}\; 
  \mho^{j_1,j_2,j}_{m_1,m_2,m}\; \Omega^{j_1,j_2,j}_{k_1,k_2,m'} 
  \Dj{j}{m'}{m}.
  \label{rel3}
\eea
\end{cor}
$Proof$ : (\ref{rel1}) has already been obtained in the proof of 
Theorem \ref{thmWigpro}, see (\ref{prod1}). Others can be 
obtained from (\ref{rel1}) by the orthogonality of $\Omega $ and $ \mho $. 
\hfill $ \Box $

  For $ {\cal G} = SL_q(2) $ and $ {\bf g} = {\cal U}_q(sl(2)) $, 
the CGC $ \Omega, \mho $ are 
given by the $q$-analogue of the CGC of $ sl(2) $ : 
$ \Omega^{j_1,j_2,j}_{m_1,m_2,m}  = \mho^{j_1,j_2,j}_{m_1,m_2,m} = 
{}_qC^{j_1,j_2,j}_{m_1,m_2,m} $ . From the relations (\ref{Wigpro}) and 
(\ref{rel1}) - (\ref{rel3}), 
the recurrence relations and the orthogonality of $ SL_q(2) $ \D are obtained 
\cite{bl,gkk,nom1}.  

  Next we show that the \D (\ref{DandT}) satisfy the RTT type relation.
\begin{thm}
The \D for $ {\cal G} $ satisfy
\beq
  \sum_{s_1,s_2}\; (R^{j_1,j_2})^{s_1,s_2}_{m_1,m_2} \Dj{j_1}{s_1}{k_1} 
  \Dj{j_2}{s_2}{k_2} 
  =
  \sum_{s_1,s_2}\; \Dj{j_2}{m_2}{s_2} \Dj{j_1}{m_1}{s_1} 
  (R^{j_1,j_2})^{k_1,k_2}_{s_1,s_2}, 
 \label{rtt}
\eeq
where $ (R^{j_1,j_2})^{s_1,s_2}_{m_1,m_2} $ are the matrix elements of the 
universal $R$-matrix for $ {\bf g} $
\[
  (R^{j_1,j_2})^{s_1,s_2}_{m_1,m_2} = \bra{j_1 m_1} \otimes \bra{j_2 m_2} 
  {\cal R} \ket{j_1 s_1} \otimes \ket{j_2 s_2}.
\]
\label{thmrtt}
\end{thm}
$ Remark $ : For $ j_1 = j_2 = 1/2 $, the matrix elements for $ {\cal R} $ are evaluated in 
the fundamental representation of $ {\bf g} $. Therefore, the relation (\ref{rtt}) 
is reduced to the defining relation of $ {\cal G} $ in FRT-formalism \cite{frt}. This implies that 
$ \Dj{\hf}{m'}{m} $ are generators of $ {\cal G}. $ 

\noindent
$ Proof $ : The relation (\ref{rtt}) can be proved by evaluating matrix elements 
of the RTT type relation for the universal $T$-matrix \cite{bcgpst}.
We define 
\[
  {\cal T}_1 = \sum\; x^{\alpha} \otimes X_{\alpha} \otimes 1, 
  \qquad
  {\cal T}_2 = \sum\; x^{\alpha} \otimes 1 \otimes X_{\alpha},
\]
then
\bea
 & & {\cal T}_1 {\cal T}_2 = \sum_{\alpha,\beta}\; x^{\alpha}x^{\beta} 
 \otimes X_{\alpha} \otimes X_{\beta} = \sum_{\alpha}\; 
 x^{\alpha} \otimes \Delta(X_{\alpha}), 
 \nn \\
 & & {\cal T}_2 {\cal T}_1 = \sum_{\alpha,\beta}\; x^{\beta}x^{\alpha} 
 \otimes X_{\alpha} \otimes X_{\beta} = \sum_{\alpha}\; 
 x^{\alpha} \otimes \Delta'(X_{\alpha}), 
 \nn
\eea
where $ \Delta' $ stands for the opposite coproduct. It follows that
\[
 {\cal T}_2 {\cal T}_1 = \sum_{\alpha}\; x^{\alpha} \otimes 
 {\cal R} \Delta(X_{\alpha}) {\cal R}^{-1},
\]
thus we obtain
\[
  (1 \otimes {\cal R}) {\cal T}_1 {\cal T}_2 = {\cal T}_2 {\cal T}_1 
  (1 \otimes {\cal R}).
\]
Evaluating the matrix elements on $ 1 \otimes \ket{j_1 k_1} \otimes \ket{j_2 k_2} $, 
the theorem is proved. \hfill $\Box $

For $ {\cal G} = SL_q(2) $, the relaiton (\ref{rtt}) was proved by Nomura 
\cite{nom1}. However, Theorem \ref{thmrtt} shows that the relation (\ref{rtt}) 
holds for any kind of deformation of $ SL(2) $. 
In Ref.\cite{nom1}, the \D for $ SL_q(2) $ are interpreted as the wave functions 
of quantum symmetric tops in noncommutative space.

\subsection*{B. Recurrence Relations and Orthogonality-like Relations}

  In this subsection, the recurrence relations and the orthogonality-like 
relations  
of the $ SL_h(2) $ \D are derived as a consequence of the theorems 
in the previous subsection. It is known that the CGC for $ {\cal U}_h(sl(2)) $ 
are given in terms of the CGC for $ sl(2) $ and the matrix elements of
the twist element $ \F$
\beq
  \Omega^{j_1,j_2,j}_{m_1,m_2,m} = \sum_{s_1,s_2}\; 
  C^{j_1,j_2,j}_{s_1,s_2,m}\; (F^{j_1,j_2})_{m_1,m_2}^{s_1,s_2},
  \label{cgcsl2}
\eeq
where $ C^{j_1,j_2,j}_{s_1,s_2,m} $ is the CGC for $ sl(2) $ and 
$ (F^{j_1,j_2})_{m_1,m_2}^{s_1,s_2} $ is given by
\[
  (F^{j_1,j_2})_{m_1,m_2}^{s_1,s_2} = \bra{j_1,m_1} \otimes \bra{j_2,m_2} 
  \F \ket{j_1,s_1} \otimes \ket{j_2,s_2}.
\]
The explicit formula for $ (F^{j_1,j_2})_{m_1,m_2}^{s_1,s_2} $ and the 
next relation are found in Ref.\cite{na1}. 
\beq
  (F^{j_1,j_2})_{-m_1,-m_2}^{-s_1,-s_2} = ((F^{-1})^{j_1,j_2})_{s_1,s_2}^{m_1,m_2}.
  \label{fandfi}
\eeq
The CGC for $ {\cal U}_h(sl(2)) $ satisfy the orthogonality relations \cite{vdj} 
because of 
\beq
 \mho^{j_1,j_2,j}_{m_1,m_2,m} = (-1)^{j_1+j_2-j} \Omega^{\ j_1,\ j_2,\ j}_{-m_1,-m_2,-m}
 = \sum_{s_1,s_2}\; C^{j_1,j_2,j}_{s_1,s_2,m}\; 
 ((F^{-1})^{j_1,j_2})^{m_1,m_2}_{s_1,s_2}.
 \label{cgcinv}
\eeq
The relation(\ref{fandfi}) and the well known property of the $ sl(2) $ CGC are 
used in the last equality. 

Note that we have known the following fact because of the remark to Theorem 
\ref{thmrtt}.
\begin{prop}
$ \Dj{\hf}{m'}{m} $ are the generators of $ SL_h(2) $
\beq
 \left(
 \begin{array}{cc}
 \Dj{\hf}{\hf}{\hf} & \Dj{\hf}{\hf}{-\hf} \\
 \Dj{\hf}{-\hf}{\hf} & \Dj{\hf}{-\hf}{-\hf}
 \end{array}
 \right)
 = 
 \left(
 \begin{array}{cc}
 x & u \\ v & y
 \end{array}
 \right)_.
 \label{dhalf}
\eeq
\label{propdhalf}
\end{prop}

 Let us consider the case that $ j_1 $ is arbitrary and $ j_2 = 1/2 $ in order to derive the recurrence 
relations for $ SL_h(2) $ $D$-functions. In this case, the $F$-coefficients 
have a simple form
\bea
 & & (F^{j_1,\hf})_{k_1,k_2}^{m_1,\hf} = \delta_{k_1,m_1} \delta_{k_2,\hf},
     \nn \\
 & & (F^{j_1,\hf})_{k_1,k_2}^{m_1,-\hf} = \delta_{k_1,m_1} 
     (\delta_{k_2,-\hf} - 2m_1h\delta_{k_2,\hf}),
     \nn \\
 & & ((F^{-1})^{j_1, \hf})_{m_1,\hf}^{n_1,n_2} = \delta_{m_1,n_1} 
     (\delta_{n_2,\hf} + 2m_1h \delta_{n_2,-\hf}),
     \nn \\
 & & ((F^{-1})^{j_1, \hf})_{m_1,-\hf}^{n_1,n_2} = 
     \delta_{m_1,n_1} \delta_{n_2,-\hf}.
     \nn
\eea
One can use Winger's product law, expressed in the form of 
(\ref{rel1}) and (\ref{rel2}), to derive the recurrence relations  
for $ \Dj{j}{m'}{m} $ which are reduced to the known recurrence 
relations of the $ SL(2) $ \D in the limit of $ h = 0. $ 
\begin{prop}
The $ SL_h(2) $ \D satisfy the following recurrence relations

\bigskip
\noindent
\begin{tabular}{ll}
{\rm (i)} & 
$
  \sqrt{j+k} \Dj{j}{k}{m} - (2k-1)h\sqrt{j-k+1} \Dj{j}{k-1}{m} 
$
\\
 &
$
  \qquad = \sqrt{j+m} \Dj{j-\hf}{k-\hf}{m-\hf} x
  +
  \sqrt{j-m} \Dj{j-\hf}{k-\hf}{m+\hf} (u - (2m+1)h x),
$
\\
{\rm (ii)} & 
$
\sqrt{j-k} \Dj{j}{k}{m} = 
  \sqrt{j+m} \Dj{j-\hf}{k+\hf}{m-\hf} v + 
  \sqrt{j-m} \Dj{j-\hf}{k+\hf}{m+\hf} (y - (2m+1)h v),
$ 
\\
{\rm (iii)} & 
$
  \sqrt{j+n} \Dj{j}{m}{n} = \sqrt{j+m} \Dj{j-\hf}{m-\hf}{n-\hf} 
  (x + (2m-1)hv) + \sqrt{j-m} \Dj{j-\hf}{m+\hf}{n-\hf} v,
$
\\
{\rm (iv)} &
$
  \sqrt{j-n} \Dj{j}{m}{n}  + \sqrt{j+n} (2n+1) h \Dj{j}{m}{n+1} 
$
\\
 & 
$
  \qquad = \sqrt{j+m} \Dj{j-\hf}{m-\hf}{n+\hf} (u + (2m-1)hy) 
     + \sqrt{j-m} \Dj{j-\hf}{m+\hf}{n+\hf} y,
$
\\
{\rm (v)} & 
$
  \sqrt{j-k+1} \Dj{j}{k}{m} + (2k-1)h \sqrt{j+k} \Dj{j}{k-1}{m}
$
\\
 &
$
 \qquad  = \sqrt{j-m+1} \Dj{j+\hf}{k-\hf}{m-\hf} x - 
  \sqrt{j+m+1} \Dj{j+\hf}{k-\hf}{m+\hf}(u -(2m+1)hx),
$
\\
{\rm (vi)} &
$
 \sqrt{j+k+1} \Dj{j}{k}{m} = -\sqrt{j-m+1} \Dj{j+\hf}{k+\hf}{m-\hf} v 
$
\\
 &
$ 
 \qquad + \sqrt{j+m+1} \Dj{j+\hf}{k+\hf}{m+\hf} (y - (2m+1)hv),
$
\\
{\rm (vii)} &
$
 \sqrt{j-n+1} \Dj{j}{m}{n} = \sqrt{j-m+1} \Dj{j+\hf}{m-\hf}{n-\hf}(x+(2m-1)hv) 
$
\\
 &   
$ 
 \qquad  - \sqrt{j+m+1} \Dj{j+\hf}{m+\hf}{n-\hf} v,
$
\\
{\rm (viii)} &
$
 \sqrt{j+n+1} \Dj{j}{m}{n} - \sqrt{j-n}(2n+1)h \Dj{j}{m}{n+1} 
$
\\
 & 
$
 \qquad = - \sqrt{j-m+1} \Dj{j+\hf}{m-\hf}{n+\hf} (u+(2m-1)hy) 
 + \sqrt{j+m+1} \Dj{j+\hf}{m+\hf}{n+\hf} y.
$
\end{tabular}
\label{proprec}
\end{prop}
$ Proof $ : Put $ j_2 = 1/2, \ j = j_1 + 1/2 $ in the relation (\ref{rel1}), 
then
\bea
 & &  \sqrt{j_1+k_1+1} \Dj{j_1+\hf}{k_1+\hf}{m} + \sqrt{j_1-k_1+1} 
 (\delta_{k_2,-\hf}-2k_1 h \delta_{k_2,\hf}) \Dj{j_1+\hf}{k_1-\hf}{m}
 \nn \\
 &=& 
 \sqrt{j_1+m+\hf} \Dj{j_1}{k_1}{m-\hf} \Dj{\hf}{k_2}{\hf} 
 + \sqrt{j_1-m+\hf} \Dj{j_1}{k_1}{m+\hf} (\Dj{\hf}{k_2}{-\hf} - 
 (2m+1)h \Dj{\hf}{k_2}{\hf}).
 \nn
\eea
Replacing $ j_1 + \hf $ and $ k_1+\hf $ with $ j$ and $k$, respectively, we 
obtain
\bea
 & & \sqrt{j+k} \delta_{k_2,\hf} \Dj{j}{k}{m} + \sqrt{j-k+1} (\delta_{k_2,-\hf}
     -(2k-1)h \delta_{k_2,\hf}) \Dj{j}{k-1}{m} 
     \nn \\
 &=&
 \sqrt{j+m} \Dj{j-\hf}{k-\hf}{m-\hf} \Dj{\hf}{k_2}{\hf} 
 +
 \sqrt{j-m} \Dj{j-\hf}{k-\hf}{m+\hf} 
 (\Dj{\hf}{k_2}{-\hf} - (2m+1)h \Dj{\hf}{k_2}{\hf}).
 \nn
\eea
The recurrence relations (i) and (ii) are obtained 
by putting $ k_2 = 1/2 $ and $ k_2 = - 1/2 $, respectively. 

  We repeat the similar computation for the relation (\ref{rel2}). 
We put $ j_2 = 1/2, \ j = j_1 + 1/2 $ in (\ref{rel2}), then rearrange 
some variables. We obtain
\bea
 & & \sqrt{j+n} \{ \delta_{n_2,\hf} + (2n -1)h \delta_{n_2,-\hf} \} 
     \Dj{j}{m}{n}
     + \sqrt{j-n+1} \delta_{n_2,-\hf} \Dj{j}{m}{n-1} 
     \nn \\
 &=& \sqrt{j+m} \Dj{j-\hf}{m-\hf}{n-\hf} \{ \Dj{\hf}{\hf}{n_2} 
     + (2m-1)h \Dj{\hf}{-\hf}{n_2} \}
     + \sqrt{j-m} \Dj{j-\hf}{m+\hf}{n-\hf} \Dj{\hf}{-\hf}{n_2}.
     \nn
\eea
The recurrence relations (iii) and (iv) correspond to 
the cases of $ n_2 = 1/2 $ and $ n_2 = -1/2 $, respectively. 

  The recurrence relations (v) - (viii) correspond to $ j_2 = 1/2, \ 
j = j_1 - 1/2 $. In this case, the relation (\ref{rel1}) yields after 
rearrangement of variables
\bea
  & & \sqrt{j-k+1} \delta_{k_2,\hf} \Dj{j}{k}{m} - \sqrt{j+k} 
  (\delta_{k_2,-\hf} - (2k-1)h\delta_{k_2,\hf}) \Dj{j}{k-1}{m} 
  \nn \\
  &=& \sqrt{j-m+1} \Dj{j+\hf}{k-\hf}{m-\hf} \Dj{\hf}{k_2}{\hf} 
  - \sqrt{j+m+1} \Dj{j+\hf}{k-\hf}{m+\hf} (\Dj{\hf}{k_2}{-\hf} 
  - (2m+1)h \Dj{\hf}{k_2}{\hf}).
  \nn
\eea
Putting $ k_2 = 1/2 $ and $ -1/2 $, we obtain the relations (v) and (vi), 
respectively. The relation (\ref{rel2}) yields
\bea
 & & \sqrt{j-n+1}(\delta_{n_2,\hf} + (2n-1)h \delta_{n_2,-\hf}) \Dj{j}{m}{n}
 - \sqrt{j+n} \delta_{n_2,-\hf} \Dj{j}{m}{n-1}
 \nn \\
 &=& \sqrt{j-m+1} \Dj{j+\hf}{m-\hf}{n-\hf} (\Dj{\hf}{\hf}{n_2} + 
 (2m-1)h \Dj{\hf}{-\hf}{n_2} ) 
 - \sqrt{j+m+1} \Dj{j+\hf}{m+\hf}{n-\hf} \Dj{\hf}{-\hf}{n_2}.
 \nn
\eea
The recurrence relations (vii) and (viii) are obtained as the cases of 
$ n_2 = 1/2 $ and $ n_2 = -1/2 $, respectively.
\hfill $ \Box$

It is possible to obtain the explicit form of  \D for some special cases 
such as $ \Dj{j}{m'}{j}, \ \Dj{j}{j}{m} $, by solving these recurrence 
relations. However, it seems to be difficult to derive  formulae for  
$ \Dj{j}{m'}{m} $ for any values of $ j,\ m' $ and $m$. We will solve 
this problem by using the tensor operator approach in \S V. 

  The orthogonality-like relations for $ \Dj{j}{m'}{m} $ can be obtained 
from the relations (\ref{rel1}) and (\ref{rel2}).
\begin{prop}
The \D for $ SL_h(2) $ $ \Dj{j}{m'}{m} $ satisfy the orthogonality-like relations 
which are reduced to the ortogonality relations of $ SL(2) $ \D in the 
limit of $ h = 0 $.
\bea
 & & \sum_{m_1, m_2}\; (-1)^{k_1-m_1} (F^{j,j})^{m1,-m1}_{m_1,m_2} \Dj{j}{k_1}{m_1}
 \Dj{j}{k_2}{m_2} = (F^{j,j})^{k_1,-k_1}_{k_1,k_2}, \label{ortho1} \\
 & & \sum_{k_1, k_2}\; (-1)^{m_1-k_1} ((F^{-1})^{j,j})^{k_1,k_2}_{k_1,-k_1} 
 \Dj{j}{k_1}{m_1} \Dj{j}{k_2}{m_2} = ((F^{-1})^{j,j})^{m_1,m_2}_{m_1,-m_1}.
 \label{ortho2}
\eea
\label{proporth}
\end{prop}
$ Proof $ : Consider the cases of $ j = 0, j_1 = j_2 $ in the relations 
(\ref{rel1}) and (\ref{rel2}). Writing $ j_1 = j_2 = j $, they yield
\bea
  & & \sum_{m_1,m_2}\; \Omega^{\; j,\;j,\;0}_{m_1,m_2,0} \Dj{j}{k_1}{m_1} \Dj{j}{k_2}{m_2} 
      = \Omega^{\; j,\;j\;,0}_{k_1,k_2,0},
  \nn \\
  & & \sum_{k_1,k_2}\; \mho^{\;j,\;j,\;0}_{k_1,k_2,0} \Dj{j}{k_1}{m_1} \Dj{j}{k_2}{m_2} 
      = \mho^{\; j,\;j,\;0}_{m_1,m_2,0}.
  \nn
\eea
The CGC are given by
\[
 \Omega^{\; j,\;j,\;0}_{m_1,m_2,0} = \sum_{s}\; C^{j,j,0}_{s,-s,0} \; 
 (F^{j,j})^{s,-s}_{m_1,m_2},\quad
 \mho^{\; j,\;j,\; 0}_{m_1,m_2,0} = \sum_{s}\; C^{j,j,0}_{s,-s,0} \; 
 ((F^{-1})^{j,j})^{m_1,m_2}_{s,-s},
\]
and 
\[
(F^{j,j})^{s,-s}_{m_1,m_2} = \delta_{s,m_1} \bra{j m_2} e^{-s \sigma} \ket{j\; -s}, \
((F^{-1})^{j,j})^{m_1,m_2}_{s,-s} = \delta_{s,m_1} \bra{j\; -s} e^{m_1 \sigma} 
\ket{j m_2}.
\] 
Then the proof of Proposition \ref{proporth} is straightforward. 
\hfill $\Box$

%
%
%
%
\setcounter{equation}{0}
\section{Review of {\boldmath $ SL(2) $} Representation Functions}

 This section is devoted to a review of the \D for Lie group $ SL(2). $ 
Especially, we focus on tensor operator properties and the relationship to 
Jacobi polynomials. We write the \D for $ SL(2) $ in terms of boson operators 
for the viewpoint of tensor operators.

 Let $ a_i^j,\ \bar a_i^j, \ i,j \in \{1, 2\}  $ be four copies of a 
boson operator commuting 
one another, $i.e.$,
\beq
  [\bar a_i^j,\; a_k^{\ell}] = \delta_{i,k}\delta^{j,\ell}, 
  \qquad
  [a_i^j,\; a_k^{\ell}] = [\bar a_i^j,\; \bar a_k^{\ell}] = 0.
  \label{bosons}
\eeq
It is known that the Lie algebra $ gl(2) \oplus gl(2) $ is realized by 
these boson operators. The left (lower) generators are defined by
\beq
  E_{ij} = a_i^1 \bar a_j^1 + a_i^2 \bar a_j^2,   \label{lowgen}
\eeq
the right (upper) generators are defined by
\beq
  E^{ij} = a_1^i \bar a_1^j + a_2^i \bar a_2^j.   \label{upgen}
\eeq
Then both left and right generators satisfy the $ gl(2) $ commutation 
relations and furthermore $ [E_{ij},\; E^{k,\ell}] = 0. $ 
Each $ gl(2) $ has  decomposition $ gl(2) = sl(2) \oplus u(1) $. 
The left and right $ sl(2) $ are generated by
\beq
  J_+ = E_{21}, \qquad J_- = E_{12}, \qquad J_0 = E_{22} - E_{11},
  \label{leftsl}
\eeq
and
\beq
 K_+ = E^{12}, \qquad K_- = E^{21}, \qquad K_0 = E^{11} - E^{22},
 \label{rightsl}
\eeq
respectively, and $ u(1) $ sectors by $ Z_L = -E_{11} - E_{22} $ and 
$ Z_R = E^{11} + E^{22} $. This choice of generators may be 
different from the usual one (see for example Ref.\cite{bl} \S4.4). 
However this is a suitable choice for 
twisting discussed in the next section. Note also that, in this realization, 
$ Z_L = - Z_R. $ Therefore, strictly speaking, this realization is not 
the direct sum of two copies of $ gl(2) $. 

  The \D for Lie group $ GL(2) $ can be given in terms of $ a_i^j $
\beq
 \Dclass{j}{m'}{m} = \{(j+m')!(j-m')!(j+m)!(j-m)!\}^{1/2} \sum_{K,L,M,N}
 \frac{(a_1^1)^K (a_2^1)^L (a_1^2)^M (a_2^2)^N}{K! L! M! N!},
 \label{D0} 
\eeq
where the sum over $ K, L, M $ and $N$ runs nonnegative integers 
provided that 
\beq
\begin{array}{ll}
 K + L = j+ m, & \qquad M + N = j-m, \\
 K + M = j+ m',& \qquad L + N = j-m.'
\end{array}
\label{condition}
\eeq
We obtain $ SL(2) $ \D by imposing $ a_1^1 a_2^2 - a_2^1 a_1^2 = 1 $. 

  It is not difficult to see that \D (\ref{D0}) form the 
irreducible tensor operators for both left and right $ gl(2) $, $i.e.$, 
\bea
 & & [J_{\pm},\; \Dclass{j}{m'}{m}] = \sqrt{(j\pm m') (j\mp m'+1)} 
 \Dclass{j}{m'\mp 1}{m}, 
 \nn \\
 & & [J_0, \; \Dclass{j}{m'}{m}] = -2m' \Dclass{j}{m'}{m}, 
 \qquad
 [Z_{L}, \; \Dclass{j}{m'}{m}]  = -2j \Dclass{j}{m'}{m},
 \label{lefttensor}
\eea
and 
\bea
 & & [K_{\pm},\; \Dclass{j}{m'}{m}] = \sqrt{(j\mp m) (j\pm m+1)} 
 \Dclass{j}{m'}{m\pm 1}, 
 \nn \\
 & & [K_0, \; \Dclass{j}{m'}{m}] = 2m \Dclass{j}{m'}{m}, 
 \qquad
 [Z_{R}, \; \Dclass{j}{m'}{m}]  = 2j \Dclass{j}{m'}{m}.
 \label{righttensor}
\eea

  It is well known that the \D for $ SL(2) $ can be expressed in 
terms of the Jacobi polynomials. The Jacobi polynomials are 
defined by
\beq
  P_n^{(\alpha,\beta)}(z) = \sum_{r \geq 0}\; 
  \frac{(-n)_r (\alpha+\beta+n+1)_r}{(1)_r (\alpha+1)_r} z^r, 
  \label{Jacobi}
\eeq
where $ (\alpha)_r $ stands for the sifted factorial
\[
 (\alpha)_r = \alpha (\alpha +1) \cdots (\alpha + r-1).
\]
For the case of $ SL(2) $, we have the relation $ a_1^1 a_2^2 = 1 + a_2^1 a_1^2 $. 
Using this, the \D are expressed for $ m'+m \geq 0, \ m' \geq m $
\beq
  \Dclass{j}{m'}{m} = \left\{
  \left(
  \begin{array}{c}
  j+m' \\ m'-m
  \end{array}
  \right) 
  \left(
  \begin{array}{c}
  j-m \\ m'-m
  \end{array}
  \right)
  \right\}^{1/2} 
  (a_1^1)^{m'+m} (a_1^2)^{m'-m} P_{j-m'}^{(m'-m,m'+m)}(z), 
  \label{DandJ}
\eeq
where $ z \equiv - a_2^1 a_1^2 $. We have the similar relations for 
other cases. 

%
%
%
%
\setcounter{equation}{0}
\section{Representation Functions for {\boldmath $ SL_h(2) $}}
\subsection*{A. Explicit Formulae for $D$-Functions}

  We saw, in the previous section,  that the \D for $ GL(2) $ form the 
irreducible tensor operators of both left and right $ gl(2) $. 
This fact leads us to the expectation that the \D for $ GL_h(2) $ also 
form the irreducible tensor operators of left and right $ {\cal U}_h(gl(2)) $. 
It is known that the tensor operators for $ {\cal U}_h(gl(2)) $ can be 
obtained from the ones for $ gl(2) $ by twisting \cite{na1,fiore}. Therefore, 
we may obtain the \D for $ GL_h(2) $ from the one for $ GL(2) $ by twisting 
twice. The irreducible tensor operators for $ {\cal U}_h(gl(2)) $ are 
defined by replacing the comutator on the left hand side of 
(\ref{lefttensor}) and (\ref{righttensor}) 
with the adjoint action. Let $ {\bf t} $ be a any tensor operator for 
$ {\cal U}_h(gl(2)) $ and $ X \in {\cal U}_h(gl(2)) $, then the adjoint 
action of $X$ on ${\bf t}$ is defined by \cite{rs}
\beq
  {\rm ad} X({\bf t}) = m(id \otimes S)(\Delta(X)({\bf t} \otimes 1)).
  \label{adjoint}
\eeq
The tensor operators $ {\bf t} $ for $ {\cal U}_h(gl(2)) $ and 
the tensor operators $ {\bf t}^{(0)} $ for $ gl(2) $ are related 
via the twist element $ \F $ by the relation \cite{fiore} 
(see also Ref. \cite{na1}) 
\beq
  {\bf t} = m(id \otimes S) (\F ({\bf t}^{(0)} \otimes 1) \F^{-1}).
  \label{tensortwist}
\eeq
Note that $ gl(2) $ and ${\cal U}_h(gl(2)) $ have the same commutation 
relations so that the realization (\ref{lowgen}) and (\ref{upgen}) is 
the realization of $ {\cal U}_h(gl(2)) $ as well. We consider the tensor 
operators under this realization of $ {\cal U}_h(gl(2)). $ 

  Let us first consider the simplest case : $ j = 1/2. $ 
What we obtain in this case from (\ref{D0}), (\ref{lefttensor}) and 
(\ref{righttensor}) is that the pairs $ (a_1^1, a_2^1),\ (a_1^2, a_2^2) $ 
are spinors of the left $ gl(2) $ and the pairs 
$ (a_1^1, a_1^2),\ (a_2^1, a_2^2) $ are spinors of the right $ gl(2). $ 
Namely, each boson operator $ a_i^j $ is a component of spinor for 
both left and right $ gl(2) $. This fact tells us that, by twisting 
via the elements
\beq
 \F_L = \exp\left(-\hf J_0 \otimes \sigma_L \right)_, 
 \qquad 
 \F_R = \exp\left(-\hf K_0 \otimes \sigma_R \right)_.   \label{LRtwist}
\eeq
with $ \sigma_L = -\ln(1-2hJ_+), \ \sigma_R = -\ln(1-2hK_+) $, we obtain a 
element of spinor for both left and right $ {\cal U}_h(sl(2)). $ To this end, 
it is convenient to rewrite the relation (\ref{tensortwist}) into different form. 
Let us write the twist element and its inverse as
\[
 \F = \sum_a \; f^a \otimes f_a, \qquad 
 \F^{-1} = \sum_a \; g^a \otimes g_a,
\]
then
\[
  \mu = \sum_a\; f^a S_0(f_a), \qquad 
  \mu^{-1} = \sum_a\; S_0(g^a) g_a.
\]
Noting the identity
\[
  \sum g^b \mu S_0(g_b) = \sum g^b f^a S_0(g_b f_a) = m(id \otimes S_0)(\F^{-1} \F) = 1,
\]
the relation (\ref{tensortwist}) yields
\beq
 {\bf t} = \sum f^a {\bf t}^{(0)} g^b S(f_a g_b) = \sum f^a {\bf t}^{(0)} \mu 
 S_0(f_a g_b) \mu^{-1} = \sum f^a {\bf t}^{(0)} S_0(f_a) \mu^{-1}.
 \label{tensortwist2}
\eeq
From (\ref{tensortwist2}), the twisting by $ \F_L $ reads
\[
  \sum_{n=0}^{\infty}\; \frac{1}{n!} \left(-\hf\right)^n J_0^n a_i^j 
  S_0(\sigma_L) \mu^{-1} 
  =
  a_i^j \sum_{k=0}^{\infty}\; \frac{(-1)^{ik}}{k!} \left(-\hf\right)^k S(\sigma_L)^k 
  = a_i^j \exp\{ (-1)^i \sigma_L/2 \}.
\]
We used the fact $ S(\sigma_L) = -\sigma_L $  in the last 
equality. To twist the above obtained result by $ \F_R $, we can 
repeat the similar computation. Then we have the doubly twisted boson operators
\beq
  a_i^j \exp\{ (-1)^i \sigma_L/2 + (-1)^{j+1} \sigma_R /2 \}. 
  \label{twistedboson}
\eeq
The commutation relations of the twisted boson operators (\ref{twistedboson}) 
are obtained by straightforward computation and it shows that the twisted boson 
operators give a realization of generators of $ GL_h(2) $. 
\begin{prop}
Let 
\beq
\begin{array}{ll}
 x = a_1^1 e^{(-\sigma_L + \sigma_R)/2}, & \qquad 
 u = a_1^2 e^{-(\sigma_L + \sigma_R)/2},
 \\
 v = a_2^1 e^{(\sigma_L + \sigma_R)/2}, & \qquad
 y = a_2^2 e^{(\sigma_L - \sigma_R)/2},
\end{array}
\label{bosonreal}
\eeq
then, $x,\; u,\; v $ and $ y$ satisfy the commutation realtions of the 
generators of $ GL_h(2) $ (\ref{SLh2}). 
In this realization, the central element $D$ is given 
\beq
 D \equiv xy -uv - hxv = a_1^1 a_2^2 - a_2^1 a_1^2.
 \label{hdeterminant}
\eeq
Note that the central element $D$ remains undeformed in this realization.
\label{propbosonreal}
\end{prop}
$ Proof $ : One can verify the commutation relations directly. We here give 
some useful commutation relations for the verification. 
The commutation relations between $ \sigma_L, \; \sigma_R $ and boson operators.
\begin{center}
\begin{tabular}{ll}
  $ [\sigma_L, \; a_1^1] = 2h e^{\sigma_L} a_2^1$, & \qquad
  $ [\sigma_L, \; a_1^2] = 2h e^{\sigma_L} a_2^2$,
  \\
  $ [\sigma_R,  a_1^2] = 2h e^{\sigma_R} a_1^1$, & \qquad
  $ [\sigma_R,  a_2^2] = 2h e^{\sigma_R} a_2^1$.
\end{tabular}
\end{center}
These are easily verified by using the power series expansion of 
$ \sigma_L, \sigma_R $ : 
$ \sigma_L = \displaystyle{ \sum_{n=1}^{\infty}\frac{(2hJ_+)^n}{n}_. } $
These relations can be used to prove the following commutation relations 
which hold for any real $k$
\begin{equation}
\begin{tabular}{ll}
 $ [e^{k\sigma_L},\; a_1^1] = 2hk e^{(k+1) \sigma_L} a_2^1,$ & \qquad 
 $ [e^{k\sigma_L},\; a_1^2] = 2hk e^{(k+1) \sigma_L} a_2^2,$
 \\
 $ [e^{k\sigma_R},\; a_1^2] = 2hk e^{(k+1) \sigma_R} a_1^1,$ & \qquad
 $ [e^{k\sigma_R},\; a_2^2] = 2hk e^{(k+1) \sigma_R} a_2^1. $
\end{tabular}
\label{sigmacom}
\end{equation}
\hfill $\Box$

  Next let us consider the twisting of $ \Dclass{j}{m'}{m} $ for any values of 
$j$ by the twist elements $ \F_L, \; \F_R. $ We denote the doubly twisted 
$ \Dclass{j}{m'}{m} $ by $ \Dj{j}{m'}{m} $, since it will be shown later that 
this $ \Dj{j}{m'}{m} $ gives the \D for $ GL_h(2) $. The computation is almost 
same as the case of spinors. What we need to compute is the twisting 
of $ (a_1^1)^K (a_2^1)^L (a_1^2)^M (a_2^2)^N $ in the expression (\ref{D0}). 
The twisting by $ \F_L $ reads
\bea
 & & \sum_{n=0}^{\infty}\; \frac{1}{n!} \left(-\hf \right)^n J_0^n 
 (a_1^1)^K (a_2^1)^L (a_1^2)^M (a_2^2)^N 
 S_0^n(\sigma_L) \mu^{-1}
 \nn \\
 &=& (a_1^1)^K (a_2^1)^L (a_1^2)^M (a_2^2)^N 
 \sum_{k=0}^{\infty}\; \frac{1}{k!} \left(-\hf \right)^k 
 (-K+L-M+N)^k \mu S_0^k(\sigma_L) \mu^{-1}
 \nn \\
 &=& (a_1^1)^K (a_2^1)^L (a_1^2)^M (a_2^2)^N \exp\{-(K-L+M-N) \sigma_L/2 \}.
 \nn
\eea
Further twisting by $ \F_R$ gives
\beq
   (a_1^1)^K (a_2^1)^L (a_1^2)^M (a_2^2)^N \exp\{ 
   -(K-L+M-N) \sigma_L/2 + (K+L-M-N) \sigma_R/2 \}.  \label{twistedD}
\eeq
Because of the condition (\ref{condition}), we have $ K-L+M-N = 2m' $ and 
$ K+L-M-N = 2m. $ Thus the exponential factor appeared in (\ref{twistedD}) is 
factored out the sum over $ K, L, M$ and $N$. Therefore we have proved the 
following proposition.
\begin{prop}
In the realization (\ref{lowgen}), (\ref{upgen}), 
the irreducible tensor operators of both left and right $ {\cal U}_h(gl(2)) $ 
are given by 
\beq
  \Dj{j}{m'}{m} = \Dclass{j}{m'}{m} \; e^{-m'\sigma_L + m \sigma_R}.
  \label{tensorUhgl2}
\eeq
\label{proptensorUhgl2}
\end{prop}

  One can write $ \Dj{j}{m'}{m} $ of Proposition \ref{proptensorUhgl2} in 
terms of the generators of $ GL_h(2) $ by making use of Proposition 
\ref{propbosonreal}. For real $A, \; B$, 
\bea
 & & (a_1^1)^K e^{(A \sigma_L + B \sigma_R)/2} = 
 (a_1^1)^{K-1} x e^{(A+1)\sigma_L /2 + (B-1)\sigma_R/2} 
 \nn \\
 &=& (a_1^1)^{K-1} e^{(A+1)\sigma_L /2 + (B-1)\sigma_R/2} 
 \{ e^{-(A+1)\sigma_L/2} x e^{(A+1)\sigma_L/2} \}. \nn
\eea
The expression $ \{ \cdots \} $ in the last line can be calculated by using 
(\ref{sigmacom}) and gives $ x - h(A+1) v$. Thus we have obtained
\bea
  & & (a_1^1)^K e^{(A \sigma_L + B \sigma_R)/2} 
  = e^{(A+K)\sigma_L /2 + (B-K)\sigma_R/2} 
  \nn \\
  & & \qquad \qquad \times (x - h (A+K)v) (x-h(A+K-1)v) \cdots
  (x-h(A+1)v).
  \label{a11}
\eea
Similar computation gives three other identities
\bea
 & & (a_2^1)^L e^{(A \sigma_L + B \sigma_R)/2} 
 = e^{(A-L)\sigma_L/2+(B-L)\sigma_R/2}\; v^L,
 \nn \\
 & & (a_1^2)^M e^{(A \sigma_L + B \sigma_R)/2} = 
 e^{(A+M)\sigma_L/2+(B+M)\sigma_R/2} 
 \nn \\
 & & \qquad \qquad \times 
 (u-h(B+M)x -h(A+M)y + h^2 (A+M)(B+M)v)
 \nn \\
 & & \qquad \qquad \times (u-h(B+M-1)x-h(A+M-1)y+h^2(A+M-1)(B+M-1)v)
 \nn \\
 & & \qquad \times \cdots \times (u-h(B+1)x- h(A+1)y + h^2(A+1)(B+1)v),
 \label{others} \\
 & & (a_2^2)^N e^{(A \sigma_L + B \sigma_R)/2} 
 = e^{(A-N)\sigma_L/2+(B+N)\sigma_R/2} 
 \nn \\
 & & \qquad \qquad \times (y-h(B+N)v) (h-h(B+N-1)v) \cdots (y-h(B+1)v). 
 \nn
\eea
The boson operators $ a_i^j $ commute one another so that the order of 
$ a_i^j$'s in $ \Dclass{j}{m'}{m} $ is irrelevant. Therefore we can have 
some different expressions of $ \Dj{j}{m'}{m} $ according to the choice of 
the order of boson operators. We here give two of them and shall show 
that they are the representation functions of $ GL_h(2) $. 
%
%
\begin{prop}
The \D for $ GL_h(2) $ are given by
\beq
 \Dj{j}{m'}{m} = \{(j+m')!(j-m')!(j+m)!(j-m)!\}^{1/2} \sum_{K,L,M,L} 
 \frac{ X_{K} \, v^L U_{K,L,M} Y_{K,L,M,N}}{K! M! L! N!},
 \label{Dfunc1}
\eeq
where $ X_K,\; U_{K,L,M} $ and $ Y_{K,L,M,N} $ are defined by
\bea
 X_K &=& x(x+hv) \cdots (x+h(K-1)v), \nn
 \\
 U_{K,L,M} &=& (u-h(K+L)x+h(K-L)y-h^2(K^2-L^2)v) 
 \nn \\
 & \times & (u-h(K+L-1)x+h(K-L+1)y-h^2(K^2-(L-1)^2)v)
 \nn \\
 &\times & \cdots 
 \nn \\
 & \times & (u-h(K+L-M+1)x+h(K-L+M-1)y \nn \\
 & & \hspace{1cm} - h^2(K^2-(L-M+1)^2)v)
 \nn \\
 & & \nn \\
 Y_{K,L,M,N} &=& (y-h(K+L-M)v) (y-h(K+L-M-1)v)
 \nn \\
 &\times& \cdots (y-h(K+L-M-N+1)v). \nn
\eea
The \D have another expression which is 
\beq
 \Dj{j}{m'}{m} = \{(j+m')!(j-m')!(j+m)!(j-m)!\}^{1/2} \sum_{K,L,M,L}\; 
 \frac{U_{M} X_{K,M} Y_{K,M,N}\; v^L}{K! M! L! N!},
 \label{Dfunc2}
\eeq
where $ U_M, X_{K,M}, Y_{K,M,N} $ are defined by
\bea 
 U_M &=& u(u+h(x+y)+h^2v) \cdots (u+h(M-1)(x+y)+h^2(M-1)^2v),
 \nn \\
 X_{K,M} &=& (x+hMv)(x+h(M+1)v) \cdots (x+h(K+M-1)v),
 \nn \\
 Y_{K,M,N} &=& (y-h(K-M)v)(y-h(K-M-1)v) \cdots (y-h(K-M-N+1)v) v^L.
 \nn 
\eea
The sum over $ K, L, M $ and $N$ runs nonnegative integers under the 
condition (\ref{condition}).
\label{propDfunc}
\end{prop}
$ Remark $ : We obtain the \D for $ SL_h(2) $ by 
putting $ D = xy - uv - h xv = 1.$

\noindent
$ Proof $ : These expressions are obtained by using (\ref{a11}) and 
(\ref{others}). The expression (\ref{Dfunc1}) corresponds to the boson 
ordering $ (a_1^1)^k (a_2^1)^L (a_1^2)^M (a_2^2)^N $, while 
the expression (\ref{Dfunc2}) corresponds to 
$ (a_1^2)^M (a_1^1)^K (a_2^2)^N (a_2^1)^L. $ 

 To show that 
$ \Dj{j}{m'}{m} $ are the representation functions of $ GL_h(2) $, 
we must verify (\ref{dfunc}). It is obvious that $ \Dj{j}{m'}{m} \in GL_h(2) $ 
and the counit of $ \Dj{j}{m'}{m} $ is easily verified by using 
$ \epsilon(x) = \epsilon(y) = 1, \ \epsilon(u) = \epsilon(v) = 0 $. However, 
it seems to be difficult to verify the coproduct of $ \Dj{j}{m'}{m} $ by 
straightforward computation. Instead of verifying the coproduct, 
we show that $ \Dj{j}{m'}{m} $ 
satisfy the recurrence relations of Proposition \ref{proprec}. 
Note that the recurrence relations of Proposition \ref{proprec} are 
for $ SL_h(2) $. The Jordanian deformation of the Lie algebra 
$ gl(2) $ considered in this paper 
is the direct sum of the deformed $ sl(2) $ and undeformed $ u(1) $ : 
$ {\cal U}_h(gl(2)) = {\cal U}_h(sl(2)) \oplus u(1). $ This implies that the 
CGC for $ {\cal U}_h(sl(2)) $ also give the CGC for $ {\cal U}_h(gl(2)) $. 
Therefore the \D for $ GL_h(2) $ also satisfy the recurrence relations of 
Proposition \ref{proprec}. 

  As an example, we show that the $ \Dj{j}{m'}{m} $ give the solutions 
to the recurrence relation (ii). We substitute the expression (\ref{Dfunc2}) 
of the \D into the first term of the RHS of the relation (ii), then 
replace the dummy index $ L $ with $ L-1 $. It follows that 
\[
 \sqrt{j+m} \Dj{j-\hf}{k-\hf}{m-\hf}v = 
 \{(j+m)!(j-m)!(j+k-1)!(j-k)!\}^{1/2} \sum_{K,L,M,N} L \;
 \frac{U_{M} X_{K,M} Y_{K,M,N}\; v^L}{K! M! L! N!}_,
\]
where the indices $ K, L, M $ and $ N$ satifsy the condition
\beq
  K+L = j+m, \ M+N = j-m,\ K+M = j+k-1,\ L+N=j-k+1.
  \label{sumcon}
\eeq
For the second term of the RHS of the relation (ii), we use the 
expression (\ref{Dfunc1}). Replacing the index $ N $ with $ N-1 $, we 
obtain
\bea
 & & \sqrt{j-m} \Dj{j-\hf}{k-\hf}{m+\hf} (u-(2m+1)hv) \nn \\
 & & \qquad = 
 \{(j+m)!(j-m)!(j+k-1)!(j-k)!\}^{1/2} \sum_{K,L,M,N} 
 N \; \frac{ X_{K} \, v^L U_{K,L,M} Y_{K,L,M,N}}{K! M! L! N!}_, \nn
\eea
where the indices $ K, L, M $ and $ N $ also satisfy the 
condition (\ref{sumcon}). Since the expressions (\ref{Dfunc1}) 
and (\ref{Dfunc2}) are the different expressions of the same 
$D$-functions, it holds that 
$ U_{M} X_{K,M} Y_{K,M,N}\; v^L = 
X_{K} \, v^L U_{K,L,M} Y_{K,L,M,N}. $ Therefore the RHS of (ii) reads
\bea
 & & \{(j+m)!(j-m)!(j+k-1)!(j-k)!\}^{1/2} \sum_{K,L,M,N} 
 (L+N) \; \frac{ X_{K} \, v^L U_{K,L,M} Y_{K,L,M,N}}{K! M! L! N!}
 \nn \\
 &=& \sqrt{j-k+1} \Dj{j}{k-1}{m}. \nn
\eea

  The four-term recurrence relation (i) is reduced to a three-term relation, 
by eliminating $ \Dj{j}{k-1}{m} $ from (i) and (ii). This recurrence relation 
is easily solved by using the another expression of $ \Dj{j}{k}{m} $ corresponding to 
another ordering of boson operators. The suitable expressions for solving 
it are the ones obtained from the ordering $ (a_1^2)^M (a_2^2)^N (a_2^1)^L 
(a_1^1)^K $ and $ (a_1^1)^K (a_2^2)^N (a_2^1)^L (a_1^2)^M. $ In this way, 
we can verify the $ \Dj{j}{m'}{m} $ obtained in this Proposition 
solve all the recurrence relations given in 
Proposition \ref{proprec}. \hfill $\Box$ 

  Both expression of (\ref{Dfunc1}) and (\ref{Dfunc2}), of course, give 
the generators of $ GL_h(2) $ for $ j = 1/2 $ which reflects Proposition 
\ref{propdhalf}. The \D for $ j=1 $ reads
\beq
 {\cal D}^1 = \left(
 \begin{array}{ccc}
 x^2 + hxv & \sqrt{2} (ux+huv) & u^2 + h(ux + uy + huv) \\
 \sqrt{2} xv & D + 2uv & \sqrt{2} (uy + huv) \\
 v^2 & \sqrt{2} yv & y^2 + hyv
 \end{array}
 \right)_.
 \label{done}
\eeq
For $ SL_h(2),\ i.e.\, $ putting $ D = 1 $, 
this coincides with the one obtained by using $h$-symplecton or quantum 
$h$-plane \cite{na1}. 
Chakrabarti and Quesne obtained the $ {\cal D}^1 $ for 
two-parametric Jordanian deformation of $ GL(2) $ in the coloured representation 
through a contraction technique to the \D for standard $(q,\lambda)$-deformation 
of $ GL(2) $  \cite{cq}. To compare the present $ {\cal D}^1 $ with 
the one given in Ref.\cite{cq}, put $ \alpha = 0, \ z = 1 $ in Eqs.(4.20) and 
(4.21) of Ref.\cite{cq}. Then we see that the $D$-functions for $j=1$ of 
Ref. \cite{cq} are different form (\ref{done}).  This difference stems from 
the different choice of the basis of $ {\cal U}_h(sl(2)) $. In Ref.\cite{cq}, 
the basis introduced by Ohn \cite{ohn} is used, that is, the commutation relations 
of the generators of $ {\cal U}_h(sl(2)) $ are not same as those of $ sl(2) $. 
While the basis of this paper satisfy the same commutation relations as 
$ sl(2) $. This results the different CGC for the same algebra so that the 
recurrence relations for the \D have the different form. The CGC for Ohn's 
basis are found in Ref.\cite{vdj}. Repeating the same procedure as \S III.B, 
we obtain another form of recurrence relations. It may be easy to 
verify that the $ {\cal D}^1 $ of Ref.\cite{cq} solves these 
recurrence relations.

\subsection*{B. {\boldmath $SL_h(2)$ } $D$-Functions and Jacobi Polynomials}

  The purpose of this subsection is to show that the \D for $ SL_h(2) $ can 
be expressed in terms of Jacobi polynomials. To this end, we return to the 
boson realization of \D (Proposition \ref{proptensorUhgl2}) and use the 
fact that the \D for Lie group $ SL(2) $ are written in terms of Jacobi 
polynomials. Recall the following two facts : (1) the central element $D$ of 
$ GL_h(2) $ is not deformed in the boson realization, 
Eq.(\ref{hdeterminant}), (2) Jacobi polynomials appeared in the \D for 
$ SL(2) $ are power series in the variable $ z = -a_2^1 a_1^2. $ 
We write the \D $ \Dclass{j}{m'}{m} $ for $ SL(2) $ appeared in 
(\ref{tensorUhgl2}) in terms of Jacobi polynomials then use the easily 
proved relation 
$ (a_2^1 a_1^2)^r = (uv)^r $ 
in order to replace the variable $z = -a_2^1 a_1^2 $ with the 
$h$-deformed one $ z = -uv $.  Let us consider, as an example, the case of 
$ m'+m \geq 0, \ m' \geq m $. The $ \Dclass{j}{m'}{m} $ are given 
by (\ref{DandJ}). We rearrange the order of $ a_1^1, \ a_1^2 $ and 
$ P_{j-m'}^{(m'-m, m'+m)}(z) $ to be 
$ P_{j-m'}^{(m'-m, m'+m)}(z) (a_1^2)^{m'-m} (a_1^1)^{m'+m}. $ 
Using the relations (\ref{a11}) and (\ref{others}), we see that 
\bea
 & & (a_1^2)^{m'-m} (a_1^1)^{m'+m} e^{-m'\sigma_L + m\sigma_R} \nn \\
 &=& u (u+h(x+y)+h^2v) \cdots (u+h(m'-m-1)(x+y)+h^2(m'-m-1)^2 v) \nn \\
 & \times & (x+h(m'-m)v) (x+h(m'-m-1)v) \cdots (x+h(2m'-1)v). \nn
\eea
This completes the expression of \D in terms of Jacobi polynomials. 

  Repeating this process for other cases, we can prove the next proposition.
%
\begin{prop}
The \D for $ SL_h(2) $ are written in terms of Jacobi polynomials 
as follows:

\noindent 
{\rm (i)} $ m'+m \geq 0, \ m' \geq m $
\bea
 \Dj{j}{m'}{m} &=& N_+  P_{j-m'}^{(m'-m, m'+m)}(z) \nn \\
 &\times & u (u+h(x+y)+h^2v) \cdots (u+h(m'-m-1)(x+y)+h^2(m'-m-1)^2 v) \nn \\
 & \times & (x+h(m'-m)v) (x+h(m'-m-1)v) \cdots (x+h(2m'-1)v). \nn
\eea

\noindent 
{\rm (ii)} $ m'+m \geq 0, \ m' \leq m $
\[
 \Dj{j}{m'}{m} = N_-  P_{j-m}^{(-m'+m, m'+m)}(z) 
 x (x+hv) \cdots (x+h(m'+m-1)v)\; v^{-m'+m} 
\]

\noindent 
{\rm (iii)} $ m'+m \leq 0, \ m' \geq m $
\bea
 \Dj{j}{m'}{m} &=& N_+ P_{j+m}^{(m'-m, -m'-m)}(z) \nn \\
 &\times & u (u+h(x+y)+h^2v) \cdots (u+h(m'-m-1)(x+y)+h^2(m'-m-1)^2 v) \nn \\
 & \times & (y-h(m-m')v)(y-h(m-m'-1)v) \cdots (y-h(2m+1)v). \nn
\eea

\noindent 
{\rm (iv)} $ m'+m \leq 0, \ m' \leq m $
\bea
 \Dj{j}{m'}{m} &=& N_-  P_{j+m'}^{(-m'+m,-m'-m)}(z) \nn \\
 & \times & v^{-m'+m}\; (y-h(m-m')v) (y-h(m-m'-1)v) \cdots (y-h(2m+1)v). \nn
\eea
The variable $ z $ is defined by $ z = -uv $ and the factors $ N_+, \ N_- $ by 
\[
 N_+ = \left\{
 \left(
 \begin{array}{c}
 j+m' \\ m'-m
 \end{array}
 \right) \left(
 \begin{array}{c}
 j-m \\ m'-m
 \end{array}
 \right) \right\}^{1/2}_,
 \quad
 N_- = \left\{
 \left(
 \begin{array}{c}
 j-m' \\ m-m'
 \end{array}
 \right) \left(
 \begin{array}{c}
 j+m \\ m-m'
 \end{array}
 \right) \right\}^{1/2}_.
\]
\label{propJacobi}
\end{prop}
$ Remark $ : The Jacobi polynomials are to the left of the generators of $ SL_h(2) $. 
To move $ P_n^{(\alpha,\beta)}(z) $ to the right, the relation
\[
 (uv)^r \exp( -m' \sigma_L + m \sigma_R ) = 
 \exp( -m' \sigma_L + m \sigma_R ) \{ uv - 2h(-m' yv + m xv) - 4h^2mm' v^2 \}^r,
\]
is used and we see that the Jacobi polynomials are changed to the power series 
in $ \zeta_{m',m} = -(u+2h(m'y-mx)-4h^2mm')v $, but the rests of the formulae 
remain unchanged.

%
%
%
%
\setcounter{equation}{0}
\section{Boson Realization of {\boldmath $GL_{h,g}(2)$} }

  It is natural to generalize the results in the previous section 
to the two-parametric Jordanian deformation of $ GL(2) $ \cite{agh}, 
since the twist element which generates the two-parametric Jordanian 
quantum algebra $ {\cal U}_{h,g}(gl(2)) $ \cite{adm,para} is known \cite{na3}. 
Unfortunately, the method in the previous sections leads us to quite complex 
calculation. As the first step to obtain the \D for two-parametric 
Jordanian quantum group $ GL_{h,g}(2) $, we here give the boson 
realization of the generators of $ GL_{h,g}(2) $. 

  The left and right twist elements are given by
\bea
 & & \F_L = \exp\left( \frac{g}{2h} \sigma_L \otimes Z_L \right) 
        \exp\left( -\frac{1}{2} J_0 \otimes \sigma_L \right), 
 \nn \\
 & & \F_R = \exp\left( \frac{g}{2h} \sigma_R \otimes Z_R \right) 
        \exp\left( -\frac{1}{2} K_0 \otimes \sigma_R \right),
 \nn
\eea
respectively. We can see that  the $ GL_{h,g}(2) $ is reduced to 
$ GL_h(2) $ when $ g = 0 $. Repeating the same procedure as (\ref{twistedboson}), 
we obtain the twisted boson operators. We can rewrite the twisted 
boson operator in terms of the generators $ GL_h(2) $. The next 
proposition can be regarded as a realization of $ GL_{h,g}(2) $ by 
generators of $ GL_h(2) $ and $ Z_L,\; Z_R$ as well.
%
%
\begin{prop}
Let 
\beq
\begin{array}{ll}
 a = x - g v Z_L, & \qquad 
 b = u - g x Z_R - g y Z_L + g^2 v Z_L Z_R,
 \\
 c = v, & \qquad
 d = y - g v Z_R,
\end{array}
\label{realgl2}
\eeq
where $ x, u, v $ and $ y $ are given by (\ref{bosonreal}). 
Then $ a, b, c $ and $ d $ satisfy the commutation relation of 
$ GL_{h,g}(2) $. 
\label{propgl2}
\end{prop}
$ Remark $ : In this realization, the quantum determinant 
$ D' = ad - bc - (h+g) ac $ for 
$ GL_{h,g}(2) $ and $ D $ for $ GL_{h}(2) $ coincide : 
$ D' = D = a_1^1 a_2^2 - a_2^1 a_1^2 $. 

\noindent
$ Proof $ : It requires a lengthy calculation, however, the proof 
is straightforward. The following commutation relations \cite{agh} are verified.
\bea
 & & [a,\; b] = -(h+g) (D'-a^2),  \qquad \quad \ \:
 [a,\; c] = -(h-g) c^2,
 \nn \\
 & & [a,\; d] = (h+g) ac - (h-g) dc, \qquad
 [b,\; c] = -(h+g) ac - (h-g)cd, 
 \label{GL2comm} \\
 & & [b,\; d] = (h-g)(D'-d^2), \qquad \qquad \:
 [c,\; d] = (h+g) c^2. 
 \nn
\eea
\hfill $\Box$

%
%
%
%
\setcounter{equation}{0}
\section{Concluding Remarks}

 In this paper, the explicit formulae of the \D for $ SL_h(2) $ (and 
$GL_h(2)$) have been obtained by using the tensor operator technique. 
We used the fact that the \D for Lie group $ GL(2) $ form irreducible 
tensor operators of $ gl(2) \oplus gl(2) $ in the realization (\ref{lowgen}), 
(\ref{upgen}). This kind of tensor operators are called double irreducible 
tensor operators in the literature. The \D for $ GL_h(2)$ were obtained 
via the construction of double irreducible tensor operators for $ 
{\cal U}_h(gl(2)) \oplus {\cal U}_h(gl(2)). $  Other examples of 
double irreducible tensor operators were considered for $q$-deformation 
\cite{queten,fiore2} and for Jordanian deformation \cite{queh}. 
Quesne constructed the $ GL_h(n) \times GL_{h'}(m) $ covariant 
bosonic and fermionic algebra which form the double irreducible tensor 
operators of $ {\cal U}_h(gl(n)) \oplus {\cal U}_{h'}(gl(m)) $  using 
the contraction method \cite{queh}. This may suggest,
in the case of $ n=m= 2 $ and $ h=h'$,  that the bosonic algebra 
of Quesne 
has a close relation to $ \Dj{\hf}{m'}{m}, \ i.e.,\ $ the generators of 
$ GL_h(2). $

 We also showed that the \D for $ SL_h(2) $ can be expressed in terms of 
Jacobi polynomials. Contrary to the $q$-deformed case where the little 
$q$-Jacobi polynomials appear in the \D for $ SU_q(2) $, the ordinary 
Jacobi polynomials are associated with the \D for $SL_h(2) $. 
It seems to be a general feature of Jordanian deformation that 
the ordinary orthogonal polynomials are associated with the representations. 
It is known that  the ordinary Gauss 
hypergeometirc functions are associated with $h$-symplecton \cite{na1}, 
while the $q$-hypergeometric functions are associated with the $q$-deformation 
of symplecton. 

 The extension of the results of this paper to the Jordanian deformation of $ SL(n) $ 
may be possible, since the explicit expressions for the twist element are 
known for the Lie algebra $ sl(n) $ \cite{klm}.

%
%
%
%

\end{document}